\documentclass{amsart}
\usepackage{caption}
\usepackage{subcaption}
\usepackage{amssymb}
\usepackage{graphicx}
\usepackage{wrapfig}
\begin{document}

\author{G\'abor. Fejes T\'oth}
\address{Alfr\'ed R\'enyi Institute of Mathematics,
Re\'altanoda u. 13-15., H-1053, Budapest, Hungary}
\email{gfejes@renyi.hu}

\author{W{\l}odzimierz. Kuperberg}
\address{Department of Mathematics \& Statistics, Auburn University, Auburn, AL36849-5310, USA}
\email{kuperwl@auburn.edu}

\title{Four classic problems}
\thanks{The English translation of the book ``Lagerungen in der Ebene,
auf der Kugel und im Raum" by L\'aszl\'o Fejes T\'oth will be
published by Springer in the book series Grundlehren der
mathematischen Wissenschaften under the title
``Lagerungen---Arrangements in the Plane, on the Sphere and
in Space". Besides detailed notes to the original text the
English edition contains eight self-contained new chapters
surveying topics related to the subject of the book but not
contained in it. This is a preprint of one of the new chapters.} 

\begin{abstract}
In this work we survey four classic problems: Borsuk's partition problem,
Tarski's plank problem, the Kneser--Poulsen problem on the
monotonicity of the union of balls under a contraction of their
centers, and the Hadwiger--Levi problem on covering convex bodies by their
smaller positively homothetic copies
\end{abstract}

\maketitle

\section{The Borsuk problem}

{\sc Borsuk} \cite{Borsuk} asked the question whether every bounded set in
$n$-dimensional space can be partitioned into $n+1$ subsets of smaller diameter.
Although Borsuk did not suggest a positive solution of the problem, for a long
while there was general belief that the answer is yes, so, the problem became
known as Borsuk's conjecture.

The truth of the conjecture for $n=2$ follows easily by the theorem of {\sc{P\'al}}
\cite{Pal20} stating that every set of unit diameter can be covered by a regular
hexagon of side length $1/\sqrt3$. Dissecting the hexagon into three pentagons
yields the sharp upper bound $\sqrt3/2$ for the diameter of the pieces. The
same upper bound was obtained by {\sc{Gale}} \cite{Gale53}, who obtained it by
dissecting a suitable truncation of an equilateral triangle of side-length
$\sqrt3$.

For sets of constant width some stronger theorems were proved. For a convex disk
$K$, let $d(K)$ denote the smallest number with the property that $K$ can be
covered by three sets of diameter $d(K)$. {\sc{Lenz}} \cite{Lenz56b} proved that
for a set $K$ of constant width 1 the inequality $\sqrt3-1\le d(K)\le\sqrt3/2$
holds, where equality is reached in the upper bound only for the circle, and in
the lower bound only for the Reuleaux triangle. Let $l(K)$  and $L(K)$ denote
the side-length of the largest equilateral triangle inscribed in $K$ and the
side-length of the smallest equilateral triangle circumscribed about $K$,
respectively. {\sc{Melzak}} \cite{Melzak63} proved that a set $K$ of constant
width 1 can be covered by three sets of diameter at most
$\min\{l(K),\sqrt{3}-l(K)\}$ and {\sc{Schopp}} \cite{Schopp77} proved that it
can be covered by three circular disks of diameter $L(K)/2$. Schopp also proved that
$\sqrt{12}-2\le L(K)\le\sqrt3$ for every set of constant width 1. In a theorem
of {\sc{Chakerian}} and {\sc{Sallee}} \cite{ChakerianSallee} the roles of disks
is changed to the opposite: They proved that every convex disk of unit diameter
can be covered by three copies of any set of constant width $0.9101$.

{\sc{Eggleston}} \cite{Eggleston55} settled the three-dimensional case of Borsuk's
conjecture by a rather complicated analytic argument. Simpler proofs were
given by {\sc{Gr\"unbaum}} \cite{Grunbaum57}, and {\sc{Heppes}} \cite{Heppes57}.
Both Gr\"unbaum and Heppes started with the result of {\sc{Gale}} {\cite{Gale53}
that every set of diameter 1 can be imbedded in a regular octahedron the distance
between whose opposite faces is 1. Then they observed that suitably
truncating the octahedron at three vertices, the resulting  polyhedron still can
contain every body of diameter 1, and can be partitioned into four sets of diameter
less than 1. For the diameter of the pieces Heppes proved the bound $0.9977\ldots$,
while with a more detailed analysis Gr\"unbaum got the bound $0.9885\ldots$. The
presently best known bound, $0.98$, is due to {\sc{Makeev}} \cite{Makeev97}, who
obtained it by proving that every convex body of diameter 1 is contained in a rhombic
dodecahedron with parallel faces at distance 1 apart. This last statement was also
proved by {\sc{Hausel}}, {\sc{Makai}} and {\sc{Sz\H{u}cs}} \cite{HauselMakaiSzucs}
and by {\sc{G.~Kuperberg}} \cite{KuperbergG}. {\sc{Katzarowa-Karanowa}}
\cite{Katzarowa-Karanowa} proved that every set of diameter $2$ in three-dimensional
space can be covered by four unit balls and she stated that her method can be used to
lower the radii of the balls to $0.999983$.

For finite sets of points in three dimensions a simple proof was given by
{\sc Heppes} and {\sc{R\'ev\'esz}} \cite{HeppesRevesz}. The proof of this
case follows by induction from the following conjecture
of V\'azsonyi (see {\sc{Erd\H{o}s}} \cite{Erdos46}): Among a set of $n\ge4$
points in $E^3$ there are at most $2n-2$ pairs realizing the diameter of the
set. Proofs of V\'azsonyi's conjecture were given by {\sc{Gr\"unbaum}}
\cite{Grunbaum56}, {\sc{Heppes}} \cite{Heppes56}, {\sc{Straszewicz}}
\cite{Straszewicz57}, and {\sc{Dol’nikov}} \cite{Dolnikov}. These proofs use
the {\it{ball polytope}} obtained by taking the intersection of the balls
centred at the points of the set with radius equal to the diameter of the set.
A proof avoiding the use of ball polytopes was given by {\sc{Swanepoel}}
\cite{Swanepoel08}.

{\sc{Rissling}} \cite{Rissling} proved that every set in the three-dimensional
hyperbolic space can be divided into four parts of smaller diameter and
the same is true for the spherical space for sets whose diameter is smaller
than $\pi/3$.

The theorem that Borsuk's conjecture holds for smooth convex bodies is generally
credited to {\sc{Hadwiger}} \cite{Hadwiger46}. Hadwiger only proved that Borsuk's
conjecture holds for a smooth body of constant width. From this, the validity of
Borsuk's conjecture for general smooth convex bodies follows, if we know that every
smooth convex body can be enclosed in a smooth body of constant width of the same
diameter. However, this was proved only later by {\sc{Falconer}} \cite{Falconer81}
and {\sc{Schulte}} \cite{Schulte81}. Direct proofs of the Borsuk conjecture for
smooth convex bodies were given by {\sc{Lenz}} \cite{Lenz56b} and {\sc{Melzak}}
\cite{Melzak67}. The condition of smoothness was weakened by {\sc{Dekster}}
\cite{Dekster93} who proved that the conjecture holds for every convex body for
which there exists a direction in which every line tangent to the body contains at
least one point of the body's boundary at which the tangent hyperplane is unique.

Consider a convex body $K$ in $E^n$ for which to any boundary point $x$ of $K$
there is a ball of radius $r$ contained in $K$ and containing the point $x$.
This means that a ball of radius $r$ can freely roll in $K$. {\sc{Hadwiger}}
\cite{Hadwiger48} proved the upper bound $d-2r\left(1-\sqrt{1-\frac{1}{n^2}}\right)$
for the diameter of the pieces in an optimal partition of such bodies into
$n+1$ parts. {\sc{Dekster}} \cite{Dekster89} proved a similar bound, namely
$\sqrt{d^2-r^2\frac{1-\sqrt{1-4/(n+3)}}{1+\sqrt{1+4/(n+3)}}}$ for odd $n$, and
$\sqrt{d^2-r^2\frac{1}{n+1}}$ for even $n$.

The Borsuk conjecture is known to be true for some further special cases:
{\sc{Rissling}} \cite{Rissling} proved it for centrally symmetric sets,
{\sc Rogers} \cite{Rogers71,Rogers81} for sets whose symmetry group contains
that of the regular simplex, and {\sc{Ko{\l}odziejczyk}} \cite{Kolodziejczyk}
for sets of revolution. The result of Ko{\l}odziejczyk was obtained
independently by {\sc{Dekster}} \cite{Dekster95}, who proved the same also
for hyperbolic and spherical spaces.

The obvious simplicia1 decomposition shows that $B^n$ can be decomposed into
$n+1$ subsets of diameter $\left(\frac{n+1}{n+2}\right)^{1/2}$, if $n$ is even,
and $\left(\frac{1}{2}+\frac{1}{2}\left(\frac{n-1}{n+3}\right)^{1/2}\right)^{1/2}$ if $n$ is odd.
It is conjectured that this is the lower bound for the diameter $d$ of the pieces.
{\sc{Hadwiger}} \cite{Hadwiger54} confirmed this for $n\le3$ and proved
the lower bound $d\ge\left(\frac{1}{2}+\frac{1}{2}\left(\frac{n-1}{2n}\right)^{1/2}\right)^{1/2}$
for $n\ge4$. {\sc{Larman}} and {\sc{Tamvakis}} \cite{LarmanTamvakis} improved
Hadwigrer's bound to $d\ge1-\frac{3}{2n}\log{n}+O(\frac{1}{n})$.

Let $b(n)$ denote the {\it $n$-th Borsuk number}, that is, the smallest
integer such that every bounded set in $n$-dimensional space can be
partitioned into $b(n)$ subsets of smaller diameter. {\sc{Lassak}}
\cite{Lassak82} established the upper bound $b(n)\le2^{n-1}+1$. Lassak's
result was improved significantly by {\sc{Schramm}} \cite{Schramm88} and
{\sc{Bourgain}} and {\sc{Lindenstrauss}} \cite{BourgainLindenstrauss} to
$$b(n)\le(\sqrt{3/2}+o(1))^n=(1.2247\ldots+o(1))^n,$$
which is the best presently known bound.

Despite the fact that some doubt in the truth of the conjecture was announced
by {\sc{Erd\H{o}s}} \cite{Erdos81}, {\sc{Larman}} \cite{Larman84}, and {\sc{Rogers}}
\cite{Rogers71}, it came as a surprise when {\sc Kahn} and {\sc Kalai} \cite{KahnKalai}
proved that the conjecture fails in all dimensions $n\ge2015$. Moreover,
they proved that $b(n)\ge(1.2)^{\sqrt n}$ for sufficiently large $n$. The best
lower bound for $b(n)$  presently known is
$$b(n)\ge{\left({\left(\frac{2}{\sqrt3}\right)}^{\sqrt2}+o(1)\right)}^{\sqrt{n}}={(1.2255\ldots+o(1))}^{\sqrt{n}},$$
due to {\sc{Ra\u{\i}gorodski\u{\i}}} \cite{Raigorodskii99}. {\sc Kahn} and
{\sc Kalai} \cite{KahnKalai} claimed without giving any details
that the Borsuk conjecture fails in dimension $n=1325$. {\sc{Weissbach}}
\cite{Weissbach00} pointed out that this statement does not follow from
the argument of Kahn and Kalai (see also {\sc{Jenrich}} \cite{Jenrich}).

The lower bound for the dimension $n$ in which the Borsuk conjecture fails
was lowered by {\sc{Nilli}} \cite{Nilli} to 946, by {\sc{Grey}} and
{\sc{Weissbach}} \cite{GreyWeissbach} to 903, by
{\sc{Ra\u{\i}gorodski\u{\i}}}, \cite{Raigorodskii97} to 561, by
{\sc{Weissbach}} \cite{Weissbach00} to 560, by
{\sc{Hinrichs}} \cite{Hinrichs02} to 323, by {\sc{Pikhurko}} \cite{Pikhurko}
to 321, by {\sc Hinrichs} and {\sc Richter} \cite{HinrichsRichter03} to 298,
and by {\sc{Bondarenko}} \cite{Bondarenko} to 65. The last step thus far was
done by {\sc{Jenrich}} and {\sc{Brouwer}} \cite{JenrichBrouwer}, who found a
64-dimensional subset of 352 points of the set constructed by Bondarenko that
cannot be divided into fewer than 71 parts of smaller diameter.

The notion of the $k$-fold Borsuk number of a set was introduced by {\sc
Hujter} and {\sc{L\'{a}ngi}} \cite{HujterLangi} as follows.  Let $S$ be a set
of diameter $d > 0$. The smallest positive integer  $m$ such that there is a
$k$-fold covering of $S$ with $m$ sets of diameters strictly smaller than
$d$, is called the {\it $k$-fold Borsuk number} of $S$.  Besides presenting a
few other results concerning this notion, Hujter and L\'{a}ngi determined the
$k$-fold Borsuk number for every bounded planar set. {\sc{L\'angi}} and
{\sc{Nasz\'odi}} \cite{LangiNaszodi17} investigated multiple Borsuk numbers
in normed spaces.

For comprehensive surveys of the topic see {\sc{Gr\"unbaum}} \cite{Grunbaum63b},
{\sc{Ra\u{\i}gorodski\u{\i}}} \cite{Raigorodskii04,Raigorodskii07,Raigorodskii08},
{\sc{Kalai}} \cite{Kalai15}, and the corresponding chapters of the books
{\sc{Boltjanski\u{\i}}}, {\sc{Martini}} and {\sc{Soltan,}}
\cite{BoltjanskiiMartiniSoltan} and {\sc{Martini, Montejano}} and {\sc{Oliveros}}
\cite{MartiniMontejanoOliveros}.

\section{Tarski's plank problem}

In \cite{Tarski} {\sc{Tarski}} raised the following problem: Is it true that
if a convex body $C$ of width $w$ is covered by parallel slabs with the widths
$w_1,\ldots,w_l$ then $w_1+\ldots+w_l\ge{w}$? If $C$ is a circle the solution
was given by {\sc{Moese}} \cite{Moese}. {\sc{Straszewicz}} \cite{Straszewicz48}
solved the problem in the plane for two strips. An affirmative answer to the
question for general convex bodies was given by {\sc Bang} \cite{Bang50,Bang51}.
Variations of Bang's proof were given by {\sc{Fenchel}} \cite{Fenchel51} and
{\sc{Bogn\'ar}} \cite{Bognar}.

The  width of a slab relative to a convex body $C$ is the width of the slab
divided by the width of $C$ in the direction perpendicular to the slab.
{\sc Bang} \cite{Bang50} asked whether the following generalization of his
theorem is true. If some slabs cover a convex body $C$ then the sum of the
widths of the slabs relative to $C$ is at least 1. For centrally symmetric
bodies Bang's question was answered in the affirmative
by {\sc{Ball}}~\cite{Ball91}. As a corollary Ball proved the following
theorem. Given a centrally symmetric convex body $C$ and $n$ hyperplanes in
$n$-dimensional Euclidean space, then there is a translate of $\frac{1}{n+1}C$
inside $C$ whose interior does not meet any of the hyperplanes.
The result is obviously sharp for every $n$ and $C$ and is a generalization
of a result by {\sc{Davenport}}~\cite{Davenport} who considered the special
case when $C$ is a cube. For non-symmetric sets $C$ Bang's problem was solved
only for coverings of $C$ by two slabs (see {\sc Bang} \cite{Bang54},
{\sc{Moser}} \cite{Moser}, {\sc{Alexander}} \cite{Alexander68}, and
{\sc{Hunter}} \cite{Hunter}).

Related to the above corollary is Conway's fried potato problem, phrased by
{\sc{Croft}}, {\sc{Falconer}}  and {\sc{Guy}}~\cite[Problem C1, p. 80]{CroftFalconerGuy}
as follows. ``In order to fry it
as expeditiously as possible Conway wishes to slice a given convex potato into
$n$ pieces by $n - 1$ successive plane cuts (just one piece being divided by
each cut) so as to minimize the greatest inradius of the pieces.'' This problem
was solved by {\sc{A.~Bezdek}} and {\sc{K.~Bezdek}}~\cite{BezdekABezdekK95}. In {\sc{A.~Bezdek}}
and {\sc{K.~Bezdek}}~\cite{BezdekABezdekK96} the problem is generalized and solved for the
case in which the role of the inradius is played by the maximum positive
coefficient of homothety of a given convex body contained in the slices.

{\sc{Ohmann}}~\cite{Ohmann} proved the following generalization of the planar
case of Bang's theorem : If a convex disk is covered by a finite family of
convex disks, then the sum of the inradii of the covering disks is at least
as large as the inradius of the covered domain. Theorem 1 in
{\sc{A.~Bezdek}}~\cite{BezdekA07} directly implies the same result. {\sc{Kadets}}
\cite{Kadets} extended this result to $n$ dimensions. Instead of the inradius
{\sc{Akopyan}} and {\sc{Karasev}} \cite{AkopyanKarasev}  measured a convex body
$B$ by the size $r_K(B)=\sup\{h\ge0: h\,K+t\subset{B}\}$ of the greatest
positively homothetic copy of a given convex body $K$ contained in it. They
proved that if in the plane $C_1,\ldots,C_k$ form a convex partition of the
convex disk $K$, then $\sum_{i=1}^k r_K(C_i)\ge1$. In higher dimensions they
proved the analogous statement for special partitions only. It should be
mentioned that not all coverings can be reduced to a partition. The question
whether an analogous result holds for coverings remains open.

{\sc{A.~Bezdek}} \cite{BezdekA03} made the following conjecture: For every convex disk
there exists an $\varepsilon>0$ such that the minimum total width of planks
needed to cover the annulus obtained by removing from the disk its
$\varepsilon$-homothetic copy contained in its interior is the same as for the
whole disk. He verified the conjecture for the case of a square with
$\varepsilon=1-\sqrt2/2$. He further supported the conjecture by proving that it
is true for every polygon whose incircle is tangent to two of its parallel sides.
It turned out that the conjecture stated in such generality is false, as
{\sc{White}} and {\sc{Wisewell}} \cite{WhiteWisewell} noticed. In fact, they
characterized all convex polygons for which Bezdek's conjecture holds
as the polygons with no minimum-width chord that meets a vertex
and divides the angle at that vertex into two acute angles.

{\sc{Zhang}} and {\sc{Ding}} \cite{ZhangDing} showed, with a very short proof, that the
equilateral triangle with an arbitrarily small hole placed anywhere in its
interior is a counterexample to Bezdek's conjecture. Also, they gave a
positive result for parallelograms. {\sc{Smurov, Bogataya}} and {\sc{Bogaty\u{i}}}
\cite{SmurovBogatayaBogatyi} proved that Bezdek's conjecture
holds for the cube of dimension $n\ge2$, even with infinitely many cubical
holes, each homothetic to the covered cube, if the hole's total edge length is
sufficiently small.

{\sc{L.~Fejes T\'oth}} \cite{FTL73c} considered the following problem: Place
$k$ great circles on a sphere so that the maximum inradius of the regions into
which the circles partition the sphere is as small as possible. He conjectured
that in the optimal arrangement the great circles dissect the sphere into the
regular tiling $\{2, 2k\}$ with congruent digonal faces. {\sc{Rosta}} \cite{Rosta}
confirmed the conjecture for $k=3$ and {\sc{Linhart}} \cite{Linhart74a}
proved it for $k=4$.

A great circle of the unit sphere and a number $r>0$ define a {\it{zone}} of width
$2r$ on the sphere, consisting of points that are at a distance at most $r$ from
the great circle. A different way to state the same problem is: Find the smallest
number $r_k$ such that the sphere can be covered by $k$ zones of width $2r_k$.
A zone of width $2r_k$ has area $\sin{r_k}$, thus $w_k\le\arcsin{1/k}$.
{\sc{Fodor, V\'{\i}gh}} and {\sc{Zarn\'ocz}} \cite{FodorVighZarnocz16b} gave an
improvement of this trivial bound. This reformulation of the problem gives
rise to the following more general conjecture: If the sphere is covered by a
finite number of zones, then the total width of the zones is at least $\pi$.
{\sc{Jiang}} and {\sc{Polyanskii}} \cite{JiangPolyanskii} gave a short, elegant
proof of this conjecture valid in all dimensions. {\sc{Ortega-Moreno}}
\cite{Ortega-Moreno} gave an alternative proof for the special case of zones
of the same width. As a corollary of their theorem Jiang and Polyanskii also
proved that if a centrally symmetric convex body on the sphere is covered by
zones of total width $w$, then it can be covered by one zone of width $w$.

Another spherical version of Tarski's plank problem, on covering the $n$-dimensional
spherical ball with a family of spherical convex bodies was considered by
{\sc{K.~Bezdek}} and {\sc{Schneider}} \cite{BezdekKSchneider}. In their theorem,
the inradius of a set used for the covering plays the role of the width of a
plank: If on the $n$-dimensional sphere the spherically convex bodies
$C_1,\ldots, C_m$ cover a spherical ball $B$ of radius $r(B)\ge\pi/2$, then
the sum of their inradii is greater than or equal to $r(B)$.

{\sc{Steinerberger}} \cite{Steinerberger18} gave lower bounds for the sum of
the $s$'th power of the areas of pairwise intersections of $n$ congruent zones
on $S^2$. His bound is asymptotically sharp for $0<s<1$. {\sc{Bezdek, Fodor,
V\'{\i}gh}} and {\sc{Zarn\'ocz}} \cite{BezdekAFodorVighZarnocz} investigated
the multiplicity of  points covered by zones. They showed that it is possible
to arrange $n$ congruent zones of suitable width on $S^{d-1}$ such that no
point belongs to more than a constant number of zones, where the constant
depends only on the dimension and the width of the zones. They also proved
that it is possible to cover $S^{d-1}$ by $n$ equal zones such that each point
of the sphere belongs to at most $c_d\ln{n}$ zones.

Concerning the Tarski plank problem and its generalizations, we refer the reader
to the book and survey by {\sc{K.~Bezdek}} \cite{BezdekK10,BezdekK14}.

\section{The Kneser-Poulsen problem}

\newfam\dsfam
\def\mathbb #1{{\fam\dsfam\tends #1}}
\font\twelveds=dsrom12
\font\tends=dsrom10
\font\eightds=dsrom8
\textfont\dsfam=\tends
\scriptfont\dsfam=\eightds

The following attractive problem was stated independently by {\sc{Poulsen}}
\cite{Poulsen} and {\sc{Kneser}}~\cite{Kneser}: If in $n$-dimensional
Euclidean space finitely many balls are rearranged so that no distance between
their centers increases, then the volume of their union does not increase.
The problem turned out to be more difficult than it appears, even in the
plane. The first result supporting the conjecture was obtained by
{\sc{Habicht}} (see {\sc{Kneser}}~\cite[p. 388]{Kneser}) and {\sc{Bollob\'{a}s}}
~\cite{Bollobas68b}, who proved it for congruent circular disks under the assumption
that the rearrangement is the result of a continuous motion during which
all distances between the disks' centers change monotonically.
{\sc{Csik\'os}}~\cite{Csikos97} and, independently, {\sc{Bern}} and {\sc{Sahai}}
\cite{BernSahai} extended this result to arbitrary circular disks. Soon
after he published this result, {\sc{Csik\'{o}s}}~\cite{Csikos98} generalized it to
balls in every dimensions.

Under the assumption of a continuous motion of the ball's centers the monotonicity
of the volume of the intersection also holds in spherical and hyperbolic space
(see {\sc{Csik\'{o}s}} \cite{Csikos01})}. However, {\sc{Csik\'{o}s}} and
{\sc{Moussong}} \cite{CsikosMoussong} showed that in the $n$-dimensional elliptical
(real projective) space the conjecture is false. Yet, in spite of the counterexample,
a configuration of $n+1$ balls reaches maximum volume of their union if the distances
between their centers become equal to $\pi/2$, the diameter of the space.

In his proof, Csik\'os represents the moving configuration of $N$ balls
of radii $r_1,\ldots r_N$, centered at $P_1(t),\ldots,P_N(t)$ in
$E^n$ $(0\le t\le1)$, by a single moving point
$P(t)=(P_1(t),\ldots,P_N(t))$ in the configuration space $R^{nN}$ and
he assumes that, while the distances between $P_i(t)$ and $P_j(t)$ do not
increase for $i, j =1,\ldots,N$, the function $P(t)$ is analytic. Then he
derives a formula expressing the derivative of the volume of the union of
the balls as a linear combination of the derivatives of the distances
between their centers with nonnegative coefficients. This yields directly that
the volume of the union does not increase. {\sc{K.~Bezdek}} and {\sc{Connelly}}
\cite{BezdekKConnelly02}, suitably modifying Lemma 1 of {\sc{Alexander}} \cite[p.~664]{Alexander85}
and using the volume formula derived by {\sc{Csik\'{o}s}} \cite{Csikos98}, succeeded in
confirming the planar case of the Kneser-Poulsen conjecture.

The conjecture that under a contraction of the centers the volume of the
intersection of a set of balls does not decrease was stated by {\sc{Gromov}}
\cite{Gromov} and {\sc{Klee}} and {\sc{Wagon}} \cite[Problem 3.1.]{KleeWagon}.
The special case of the conjecture for congruent circular disks and continuous
motion was established by {\sc{Capoyleas}}~\cite{Capoyleas} before {\sc{K.~Bezdek}}
and {\sc{Connelly}} \cite{BezdekKConnelly02} proved it in full generality for the plane.

For a compact set $M$ in a space of constant curvature, consider those balls $B\subset{M}$
(of possibly zero radius) that are not a proper subset of any ball $B'\subset{M}$. The set
of centers of these balls is called the {\it{center set of $M$}}. {\sc{Gorbovickis}}
\cite{Gorbovickis18} invented a method dealing with Kneser--Poulsen-type problems based on
the investigation of the properties of central sets. He proved that if on a plane of
constant curvature the union of a finite set of (not necessarily congruent) closed circular
disks has a simply connected interior, then the area of the union of these disks cannot
increase after any contractive rearrangement. We emphasise the following corollary of the
spherical case of this theorem: (i) If a finite set of circular disks on the sphere with
radii not smaller than $\pi/2$ is rearranged so that the distance between each pair of
centers does not increase, then the area of the union of the disks does not increase.
(ii) If a finite set of disks on the sphere with radii not greater than $\pi/2$ is
rearranged so that the distance between each pair of centers does not increase, then
the area of the intersection of the disks does not decrease. For two dimensions, this
generalizes the result of {\sc{K.~Bezdek}} and {\sc{Connelly}}~\cite{BezdekKConnelly04}
in which they proved the analogous statement for $n$ dimensions but only for hemispheres.
The limiting case of the application of claim (ii) when the radii approach zero yields an
alternative proof of the corresponding theorem in the Euclidean plane. With a suitable
adaptation of Gorbovickis' method {\sc{Csik\'os}} and {\sc{Horv\'ath}}
\cite{CsikosHorvath18} also proved the monotonicity of the area of the intersection
of circular disks on the hyperbolic plane.

{\sc Alexander}~\cite{Alexander85} conjectured that, under an arbitrary
contraction of the centers of finitely many congruent circles, the perimeter of
the intersection of the circles does not decrease. The proof of this conjecture
does not seem to lie within reach. {\sc{K. Bezdek, Connelly}} and {\sc{Csik\'os}}
\cite{BezdekKConnellyCsikos} settled some special cases of the conjecture, among
other cases, they proved it for four circles. The weaker result concerning
the perimeter of the convex hull of the circles was proved by {\sc{Sudakov}}
\cite{Sudakov}, rediscovered by {\sc Alexander}~\cite{Alexander85} and extended
to the hyperbolic plane and to the hemisphere by {\sc{Csik\'os}} and
{\sc Horv\'ath} \cite{CsikosHorvath18}. For the Euclidean case, a simpler proof
was given by {\sc{Capoyleas}} and {\sc{Pach}}~\cite{CapoyleasPach}, who also
established a similar result in the case where the Euclidean norm is replaced
by the maximum norm.

In dimensions higher than 2, for non-continuous contractions, there are only a
few partial results concerning the Kneser--Poulsen conjecture. In $n$ dimensions,
the conjecture was verified for $n+1$ balls by {\sc{Gromov}} \cite{Gromov}.
{\sc{K.~Bezdek}} and {\sc{Connelly}}~\cite{BezdekKConnelly02} extended
Grovov's result to at most $n+3$ balls. The conjecture is proved for special
arrangements of balls, {\it{e.g.}} for a small number of intersections or large
radii by {\sc{Gorbovickis}} \cite{Gorbovickis13, Gorbovickis14}.
{\sc{K. Bezdek}} and {\sc{Nasz\'odi}} \cite{BezdekKNaszodi18a} proved that
none of the intrinsic volumes of the intersection of $k\ge(1+\sqrt{2})^d$
congruent balls in $E^n$ decreases under {\it{uniform contractions}} of the centers,
that is under contractions where all the pairwise distances in the first set of
points are larger than all the pairwise distances in the second set of points.
{\sc{K. Bezdek}} \cite{BezdekK19} gave an alternative proof of a slightly
stronger theorem. Moreover, in \cite{BezdekK20} he proved that in a Minkowski
space, under uniform contraction of the centers, the volume of both the
intersection and the union of the balls changes monotonously.

The Kneser-Poulsen conjecture makes sense in any connected Riemannian
manifold. However, {\sc{Csik\'{o}s}} and {\sc{Kunszenti-Kov\'{a}cs}}
\cite{CsikosKunszenti-Kovacs} proved that if the conjecture is true in a
Riemannian manifold, then the manifold must be of constant curvature. Moreover,
{\sc{Csik\'os}} and {\sc{Horv\'ath}} \cite{CsikosHorvath14} showed that
the same consequence holds even for balls of the same radii.

Further results and generalizations on this topic are found in {\sc{Csik\'{o}s}}
\cite{Csikos01} and {\sc{K.~Bezdek}} and {\sc{Connelly}} \cite{BezdekKConnelly08}.
The articles by {\sc{K.~Bezdek}} \cite{BezdekK08} and {\sc{Csik\'os}} \cite{Csikos18}
contain surveys on the subject of the Kneser-Poulsen conjecture.

\section{Covering a convex body by smaller homothetic copies}

{\sc{Hadwiger}} \cite{Hadwiger57d}, {\sc{Levi}} \cite{Levi55} and {\sc{Gohberg}}
and {\sc{Markus}} \cite{GohbergMarkus60}, independently of each other, asked
for the smallest integer $h(K)$ such that a given $n$-dimensional convex body
$K$ can be covered by $h(K)$ smaller positively homothetic copies of $K$.  A
boundary point $x$ of the convex body $K$ is {\it illuminated} from the
direction of a unit vector $u$ if the ray issuing from $x$ in the direction of
$u$ intersects the interior of $K$. Let $i(K)$ denote the minimum number of
directions from which the boundary of $K$ can be illuminated. The problem of
finding the maximum value of $i(K)$ was raised by {\sc{Boltjanski\u{\i}}}
\cite{Boltjanskii60} and in a slightly different, but for
compact sets equivalent form, by {\sc{Hadwiger}} \cite{Hadwiger60}.
Boltjanski\u{\i} observed that for convex bodies $h(K)=i(K)$. Both
Boltjanski\u{\i} and Hadwiger conjectured that
$$i(K)\le 2^n$$
for every $n$-dimensional convex body $K$ and that equality holds only
for parallelotopes.

Both {\sc{Levi}} \cite{Levi55} and {\sc{Gohberg}} and {\sc{Markus}}
\cite{GohbergMarkus60} verified the conjecture for the plane. {\sc{Lassak}}
\cite{Lassak86} proved the sharp result that every convex disk
can be covered by four homothetic copies with ratio $\sqrt2/2$. An extreme
example is the circle. The conjecture remains open for $n\ge 3$. In three
dimensions it was proved for centrally symmetric convex bodies by {\sc{Lassak}}
\cite{Lassak84}, for convex bodies symmetric in a plane by {\sc{Dekster}}
\cite{Dekster00b}, and for convex polyhedra with an arbitrary affine symmetry
by {\sc{K.~Bezdek}} \cite{BezdekK91}.

There is a great variety of results confirming the conjecture for special
classes of bodies in $E^n$ by establishing upper bounds for $h(K)$ or $i(K)$
smaller than $2^n$. {\sc{Schramm}} \cite{Schramm88} proved that if $K$ is a set
of constant width then $i(K)<5n\sqrt{n}(4+\ln{n})\left(\frac{3}{2}\right)^{n/2}$.
{\sc{K.~Bezdek}} \cite{BezdekK11,BezdekK12b} extended Schramm's inequality to
a wider class of convex bodies, namely for those convex bodies $K$ that are the
intersection of congruent balls with centers in $K$. {\sc{Martini}} \cite{Martini87}
proved the bound $h(K)\le3\cdot2^{n-2}$ for every zonotope other than a parallelotope.
The bound $h(K)\le3\cdot2^{n-2}$, $K$ not a parallelotope, was verified also for
zonoids by {\sc{Boltjanski\u{\i}}} and {\sc{P. S.~Soltan}} \cite{BoltjanskiiSoltan},
and {\sc{Boltjanski\u{\i}}} \cite{Boltjanskii95} extended it to an even wider
class of convex bodies, the so-called belt bodies. {\sc{Boltjanski\u{\i}}} and
{\sc{Martini}} \cite{BoltjanskiiMartini87} characterized those belt bodies $K$ for
which $h(K)=3\cdot2^{n-2}$. {\sc{K.~Bezdek}} and {\sc{Bisztriczky}}
\cite{BezdekKBisztriczky} proved the conjecture for dual cyclic polytopes.
For the dual $P$ of an $n$-dimensional cyclic polytope {\sc{Talata}}
\cite{Talata99b} proved the inequality $(n+1)(n+3)/4\le h(P)\le(n+1)^2/2$ and,
for the case when $n$ is even, he gave the sharp bound $h(P)\le(n/2+1)^2$.
{\sc{Tikhomirov}} \cite{Tikhomirov} considered convex bodies $K$ in $E^n$ with the
property that for any point $(x_1,\ldots,x_n)\in{K}$, any choice of signs
$\varepsilon_1,\ldots,\varepsilon_n\in\{-1,1\}$ and any permutation
$\sigma$ on $n$ elements
$\varepsilon_1x_{{\sigma}_1},\ldots,\varepsilon_1x_{{\sigma}_n}\in{K}$.
He proved that for sufficiently large $n$, we have $i(K)<2^n$ for every
such convex body $K$ different from a cube.
Smooth convex bodies in $E^n$ can be illuminated by $n+1$ directions,
see {\it{e.g.}} {\sc{Boltjanski\u{\i}}} and {\sc{Gohberg}}
\cite[Theorem~9, p. 61]{BoltjanskiiGohberg65}. {\sc{Dekster}}
\cite{Dekster94} extended this result by replacing the assumption
of smoothness of the body by a certain smoothness of just a single belt
of the body.

{\sc{Livshyts}} and {\sc{Tikhomirov}} \cite{LivshytsTikhomirov20a} proved that
the cube represents a strict local maximum for these problems: If a convex body,
that is not a parallelotope, is close to the cube in the Banach--Mazur metric,
then it can be covered by $2^n-1$ smaller homothetic copies of itself, and
$2^n-1$ light sources suffice to illuminate its boundary.

A simple consequence of the upper bound of Rogers for the translational
covering density combined with the Rogers--Shephard inequality for the
volume of the difference body is that
$h(K)\le\binom{2n}{n}(n\log{n}+\log{n}+5n)$
for every convex body and
$h(K)\le2^n(n\log{n}+\log{n}+5n)$
for centrally symmetric convex bodies in $E^n$ (see {\sc{Rogers}} and
{\sc{Zong}} \cite{RogersZong}). An improvement by a sub-exponential factor
was recently obtained by {\sc{Huang, Slomka, Tkocz}} and {\sc{Vritsiou}}
\cite{HuangSlomkaTkoczVritsiou}. For general convex bodies their bound is
of the order of $\binom{2n}{n}e^{-c\sqrt{n}}$ for some universal constant
$c>0$. For low dimensions better bounds were given by {\sc{Lassak}}
\cite{Lassak88} and {\sc{Prymak}} and {\sc{Shepelska}} \cite{PrymakShepelska}.
In three dimensions {\sc{Lassak}} \cite{Lassak98} proved $h(K)\le20$
which was improved by {\sc{Papadoperakis}} \cite{Papadoperakis} to $h(K)\le16$
and subsequently by {\sc{Prymak}} \cite{Prymak} to 14.

A {\it{cap body}} of a ball is the convex hull of a closed ball $B$ and
a countable set $\{p_i\}$ of points outside the ball such that for any
pair of distinct points $p_i, p_j$  the line segment $p_ip_j$
intersects $B$. The difficulty of the illumination problem is exposed
by the following example of {\sc{Nasz\'odi}} \cite{Naszodi16b}. Clearly,
we have $i(B^d)=d+1$. On the other hand, for any $\varepsilon>0$ there
is a centrally symmetric cap body $K$ of $B^d$ and a positive constant
$c=c(\varepsilon)$ such that $K$ is $\varepsilon$ close to $B^d$ and
$i(K)\ge c^d$. {\sc{Ivanov}} and {\sc{Strachan}} \cite{IvanovStrachan}
studied the illumination number of cap bodies in 3 and 4 dimensions
and proved that $i(K)\le6$ for centrally symmetric cap bodies of a ball
in $E^3$, and $i(K)\le8$ for unconditionally symmetric cap bodies of a
ball in $E^4$.

{\it{Weighted illumination}} was introduced by {\sc{Nasz\'odi}} \cite{Naszodi09}
and {\it{weighted covering}} was introduced by {\sc{Artstein-Avidan}} and
{\sc{Raz}} \cite{Artstein-AvidanRaz} and {\sc{Artstein-Avidan}} and
{\sc{Slomka}} \cite{Artstein-AvidanSlomka}. A collection of weighted light
sources illuminates a convex body $K$ if for every boundary point $x$ of $K$
the total weight of light sources illuminating $x$ is at least $1$. Similarly,
a collection of weighted bodies covers $K$ if for every point $x$ of $K$
the total weight of bodies containing $x$ is at least $1$. The {\it{weighted}}
or {\it{fractional illumination number}} $i^*(K)$ of $K$ is the infimum of the
total weight of light sources illuminating $K$. The {\it{weighted covering
number}} $h^*(K)$ of $K$ is the infimum of the total weight of smaller weighted
homothetic copies of $K$ covering $K$. {\sc{Nasz\'odi}} \cite{Naszodi09}
conjectured that $i^*(K)\le2^n$ for every $n$-dimensional convex body $K$.
He proved the conjecture for centrally symmetric bodies and established the
inequality $i^*(K)\le\binom{2n}{n}$ for general convex bodies.
{\sc{Artstein-Avidan}} and {\sc{Slomka}} \cite{Artstein-AvidanSlomka}
proved that $h^*(K)\le2^n$ for centrally symmetric convex bodies
$K\subset E^n$ with equality only for a parallelotope and
$h^*(K)\le\binom{2n}{n}$ for general convex bodies. They pointed out that
the proof of the equivalence between the illumination problem and the
Levi-Hadwiger covering problem carries over to the weighted
setting. This way they gave an alternative proof of Nasz\'odi's result.

{\sc{Bezdek}} and {\sc{L\'angi}} \cite{BezdekKLangi20a} considered the
illumination problem on the sphere. A boundary point $q$ of a convex
body $K$ in $S^d$ is illuminated from a point $p\in S^n\setminus K$
if it is not antipodal to $p$, the spherical segment with endpoints
$p$ and $q$ does not intersect the interior of $K$, but the great-circle
through $p$ and $q$ does. The {\it{illumination number}} of $K$ is the
smallest cardinality of a set that illuminates each boundary point of $K$
and lies on an $(n-1)$-dimensional great-sphere of $S^n$ which is
disjoint from $K$. Bezdek and L\'angi proved that the illumination number
of every convex polytope in $S^n$ is $n+1$ and raised the question whether
there is a convex body in $S^n$ whose illumination number is greater than
$n+1$.

The paper by {\sc{Nasz\'odi}} \cite{Naszodi18} surveys different problems
about covering, among others the Hadwiger-Levi problem.
For literature and further results concerning the illumination problem we
refer to the surveys by {\sc  K.~Bezdek} \cite{BezdekK93,BezdekK06a},
{\sc{K. Bezdek}} and {\sc{Khan}} \cite{BezdekKKhan18a},
{\sc{Boltjanski\u{\i}}} and {\sc{Gohberg}} \cite{BoltjanskiiGohberg95},
{\sc{Szab\'o}} \cite{Szabo97}, {\sc{Martini}} and
{\sc{Soltan}} \cite{MartiniSoltan} and to the book by {\sc{Boltjanski\u{\i}}},
{\sc{Martini}} and {\sc{ P. S.~Soltan}} \cite{BoltjanskiiMartiniSoltan}.
\eject

\small{
\bibliography{pack}}